\documentclass[11pt,a4paper]{amsart}
\usepackage{amsthm,amsmath}
\usepackage[all]{xy}
\title{On the Weak Lefschetz Property for Hilbert functions of almost complete intersections} 
\author{Alfio Ragusa
	\and  Giuseppe Zappal\`a}

\subjclass[2000]{13 D 40, 13 H 10}
\keywords{Almost complete intersections, Gorenstein rings, Liaison theory , Hilbert function}

\DeclareSymbolFont{rsfscript}{OMS}{rsfs}{m}{n}
\DeclareSymbolFontAlphabet{\mathrsfs}{rsfscript}

\DeclareSymbolFont{AMSb}{U}{msb}{m}{n}
\DeclareSymbolFontAlphabet{\mathbb}{AMSb}

\DeclareSymbolFont{eufrak}{U}{euf}{m}{n}
\DeclareSymbolFontAlphabet{\gothic}{eufrak}

\def\ac{\`}

\newcommand{\f}{\footnotesize}

\newcommand{\pp}{\mathbb P}

\newcommand\depth{\operatorname{depth}}

\newcommand\mci{\operatorname{mci}}

\newcommand\Aci{\operatorname{Aci}}

\newcommand\Gor{\operatorname{Gor}}
\newcommand\Regor{\operatorname{Regor}}
\newcommand\zz{{\mathbb Z}}
\newcommand{\tgz}{t_{G|Z}}

\newtheorem{thm}{Theorem}[section]
\newtheorem{lem}[thm]{Lemma}
\newtheorem{prp}[thm]{Proposition}

\theoremstyle{definition}
\newtheorem{ntn}[thm]{Notation}
\newtheorem{dfn}[thm]{Definition}

\theoremstyle{remark}
\newtheorem{rem}[thm]{Remark}

\newcounter{num}
\newcommand{\dom}{\addtocounter{num}{1}\thenum)}

\begin{document}


\begin{abstract}
It is known that all complete intersection Artinian standard graded algebras of codimension $3$ have the Weak Lefschetz Property. Unfortunately, this property does not continue to be true when you increase the number of minimal generators for the ideal defining the algebra. For instance, it is not more valid for almost complete intersection  Artinian standard graded algebras of codimension $3.$ On the other hand, the Hilbert functions of all Weak Lefschetz  Artinian graded algebras are unimodal and with the positive part of their first difference forming an $O$-sequence (i.e. are Weak Lefschetz sequences). In this paper we show that all the Hilbert functions of the almost complete intersection Artinian standard graded algebras of codimension $3$ are Weak Lefschetz sequences.
\end{abstract}



\maketitle

\section*{Introduction}
\markboth{\it Introduction}{\it Introduction}
In the last few years the Weak Lefschetz Property for Artinian standard graded algebra has been investigated very deeply. It is trivial to see that all Artinian standard graded algebra of codimension $2$ have the Weak Lefschetz Property, but already in codimension $3$ many questions are still open. One of the most import result on this direction is the proof of the Weak Lefschetz Property for all Artinian complete intersections of codimension $3$, due to T.~Harima, J.C.~Migliore, U.~Nagel, J.~Watanabe, in \cite{HMNW}. Other results on Weak Lefschetz algebras and the corresponding Hilbert functions can be found, for instance, in \cite{Mi}, \cite{MZ}. Nevertheless, again in codimension $3,$ it is enough to consider almost complete intersection ideals,  i.e. perfect ideals of height $3$ and minimally generated by $4$ elements, to produce examples of Artinian standard graded algebras not enjoying the Weak Lefschetz Property (see \cite{BK}). On the other hand there exists a characterization of the sequences which occur as Hilbert functions of Weak Lefschetz algebras. These sequences are exactly those which are unimodal and such that the positive part of their first differences is an $O$-sequence (these will be named Weak Lefschetz sequences).
The aim of this paper is to show that the Hilbert functions of all almost complete intersection Artinian standard graded algebras of codimension $3$ are Weak Lefschetz sequences, see Thorem \ref{MT}.

The paper is structured in this way. Section $1$ contains notation, terminology and basic facts about almost complete intersections and their relationships with Gorenstein ideals and regular sequences which are contained in them. The main result is placed in section $2,$ which contains also all the preparatory material in order to prove it. All this is done using a claim whose technical proof is demanded to section $3.$

\section{Notation and preliminaries} 
\markboth{\it Notation and preliminaries}
{\it Notation and preliminaries}
Let $k$ be an algebraically closed field and let $R:=k[x_1,x_2,x_3].$ We consider on $R$ the standard grading and we consider in it just homogeneous ideals, which we call simply ideals. \par 
An ideal $I_Q\subset R$ is said to be an {\em almost complete intersection} ideal of codimension $3$ if $I_Q$ is minimally generated by $4$ forms and $R/I_Q$ is an Artinian $k$-algebra. 
\par
Every almost complete intersection ideal $I_Q$ of codimension $3$ is directly linked in a complete intersection to a Gorenstein ideal $I_G\subset R.$ Indeed, if $I_Z\subseteq I_Q$ is generated by $3$ minimal generators of $I_Q,$ which form a regular sequence, then $I_G:=I_Z:I_Q$ is a Gorenstein ideal. By liaison theory (see \cite{PS} for a complete discussion on this argument) we have also $I_Q=I_Z:I_G.$ 
\par
Let $e_1\le e_2\le e_3\le e_4$ be positive integers. We define
\begin{multline*}
 \Aci(e_1,e_2,e_3,e_4):=\\ \{I_Q\mid I_Q\subset R \text{ is an almost complete intersection homogeneous ideal}\\
 \text{of codimension $3,$ minimally generated in degrees $e_1,e_2,e_3,e_4$}\}.
\end{multline*}
Let $d_1\le\ldots\le d_{2m+1}$ ($m\ge 1$) be positive integers. We define
\begin{multline*}
 \Gor(d_1,\ldots,d_{2m+1}):=\\ \{I_G\mid I_G\subset R \text{ is a Gorenstein homogeneous ideal}\\
 \text{of codimension $3,$ minimally generated in degrees $d_1,\ldots,d_{2m+1}$}\}.
\end{multline*}
Of course, if $I_Z\in\Gor(a_1,a_2,a_3)$ then $I_Z$ is a complete intersection ideal.
\par
Let $A_W:=R/I_W$ be an Artinian standard $k$-algebra. Throughout the paper we will denote by $H_{A_W}$ or, more simply, by $H_{W}$ its Hilbert function. Moreover we set
 $$\vartheta_W:=\max\{i\in\zz\mid H_{W}(i)>0\}+3$$
and
 $$\lambda_W:=\max\{i\in\zz \mid\Delta H_W(i)>0\}.$$
If $A_G:=R/I_G$ is a Gorenstein Artinian standard $k$-algebra whose Hilbert function is $H_G,$ then 
 $$H_G(i)=H_G(\vartheta_G-3-i),\text{ for every }i\in\zz,$$
and consequently
 $$\Delta H_G(i)=-\Delta H_G(\vartheta_G-2-i),\text{ for every }i\in\zz.$$
In the next remark we collect some simple facts about $\vartheta_G$ and $\lambda_G$ when $I_G$ is a Gorenstein Artinian ideal.
\begin{rem}\label{rmgz}
Let $I_G\subset k[x_1,x_2,x_3]$ be a Gorenstein Artinian ideal. Then
\setcounter{num}{0}
\begin{itemize}
	\item[\dom]$\Delta H_G(i)>0$ $\iff$ $0\le i\le\lambda_G;$ 
	\item[\dom]$\Delta H_G(i)<0$ $\iff$ $\vartheta_G-2-\lambda_G\le i\le\vartheta_G-2;$
	\item[\dom]$\vartheta_G\ge 2\lambda_G+3;$
	\item[\dom]if $\vartheta_G\ge 2\lambda_G+4$ then $\Delta H_G(i)=0$ for $\lambda_G+1\le i\le\vartheta_G-3-\lambda_G;$
	\item[\dom]if $I_Z\in\Gor(a_1,a_2,a_3)$ then
	$$\lambda_Z=\left\{
\begin{array}{cl}
	\left\lfloor\frac{\vartheta_Z-3}{2}\right\rfloor & \text{ if } a_2\le a_3\le a_1+a_2 \\
	a_1+a_2-2 & \text{ if } a_3\ge a_1+a_2-1
\end{array}\right.;$$
in particular $\lambda_Z\le a_1+a_2-2;$
  \item[\dom]if $I_Z\in\Gor(a_1,a_2,a_3)$ then
  $$\Delta H_Z(\lambda_Z)=\left\{
\begin{array}{cl}
	1 & \text{ if } a_2\le a_3\le a_1+a_2-2\text{ and }\vartheta_Z\text{ is odd}\\
	2 & \text{ if } a_2\le a_3\le a_1+a_2-2\text{ and }\vartheta_Z\text{ is even} \\
	1 & \text{ if }  a_3\ge a_1+a_2-1
\end{array}\right..$$
\end{itemize}
\end{rem}
It is well known that $\Gor(d_1,\ldots,d_{2m+1})\ne\emptyset$ iff the following Gaeta conditions hold (see \cite{Ga} for the general result and \cite{Di} for the Gorenstein version)
\begin{itemize}
	\item[1)]$\vartheta:=\sum_{i=0}^{2m+1}d_i/m$ is an integer;
	\item[2)]$\vartheta>d_i+d_{2m+3-i}$ for $2\le i\le m+1.$
\end{itemize}
If $I_G\in\Gor(d_1,\ldots,d_{2m+1})$ then $\vartheta_G=\sum_{i=0}^{2m+1}d_i/m$ and   
$I_G$ admits a graded minimal free resolution of the following type
 $$0\to R(-\vartheta_G)\to\bigoplus_{i=1}^{2m+1}R(d_i-\vartheta_G)\to\bigoplus_{i=1}^{2m+1}R(-d_i)$$
(for more on this see the beautiful structure theorem of  $3$-codimensional Gorenstein ideals due to D. Buchsbaum and D. Eisenbud in \cite{BE}). \par
Let $\delta=(d_1,\ldots,d_{2m+1}),$ such that $\Gor\delta\ne \emptyset$ and let $I_G\in\Gor\delta;$ we set
 $$B_{\delta}:=\{3\le i\le m+1\mid\vartheta_G\le d_i+d_{2m+4-i}\},$$ 
 $$C_{\delta}:=\{4\le i\le m+2\mid\vartheta_G\le d_i+d_{2m+5-i}\}$$
and 
 $$\mci\delta=\left\{\begin{array}{cl}(d_1,d_2,d_3)&\text{ if }B_{\delta}=C_{\delta}=\emptyset\\
 (d_1,d_2,d_{\max C_{\delta}})&\text{ if }B_{\delta}=\emptyset\text{ and }C_{\delta}\ne\emptyset\\
 (d_1,d_{\max B_{\delta}},d_{2m+4-\min B_{\delta}})&\text{ if }B_{\delta}\ne\emptyset
 \end{array}\right.$$

Let $\alpha=(a_1,a_2,a_3),$ $a_1\le a_2\le a_3$ and $\delta=(d_1,\ldots,d_{2m+1}),$ $d_1\le\ldots\le d_{2m+1},$
$\Gor\delta\ne\emptyset;$
we set
 $$\Regor(\alpha,\delta)=\{(I_Z,I_G)\in\Gor\alpha\times\Gor\delta\mid I_Z\subseteq I_G\}.$$ 
By~\cite{RZ2}, Theorem 3.6, we have that
 $$\Regor(\alpha,\delta)\ne\emptyset\iff\alpha\ge\mci\delta.$$
Let $\varphi:\zz \to\zz$ be a function. We set $$\varphi^+:=\frac{|\varphi|+\varphi}{2}.$$
Now, let us suppose that $\varphi(i)=0$ for every $i\in \zz^-;$ we remind that $\varphi$ is called {\em unimodal} if there exists $u\in\zz^+$ such that
 $$\begin{array}{ll}\varphi(i)<\varphi(i+1)&\text{ for }0\le i\le u-1\\
 \varphi(i)\ge\varphi(i+1)&\text{ for }i\ge u\,.\end{array}$$
\begin{prp}\label{crescdelta}
Let $\delta=(d_1,\ldots,d_{2m+1}),$ with $d_1\le\ldots\le d_{2m+1}$ ($m\ge 2$) positive integers, such that 
$\Gor\delta\ne\emptyset.$ Let $\vartheta=\sum_{i=0}^{2m+1}d_i/m.$ Let us suppose that there exist $h$ and $k,$ $h<k,$ such that $d_h+d_k=\vartheta.$ Let    
 $$\widehat{\delta}:=(d_1,\ldots,d_{h-1},d_{h+1},\ldots,d_{k-1},d_{k+1},\ldots,d_{2m+1})$$
Then
\begin{itemize}
	\item[1)]$\Gor\widehat{\delta}\ne\emptyset;$
	\item[2)]$\mci\delta\ge\mci\widehat{\delta}.$
\end{itemize}
\end{prp}
\begin{proof}
\begin{itemize}
	\item[1)]$\widehat{\delta}$ satisfies Gaeta conditions too.
	\item[2)]This is true by Theorem 3.9 in \cite{RZ2}.
\end{itemize}
\end{proof}
\begin{dfn}
Let $H$ be the Hilbert function of an Artinian standard graded $k$-algebra. We will say that $H$ is a {\em Weak Lefschetz sequence} if it is unimodal and $(\Delta H)^+$ is an $O$-sequence.
\end{dfn}
For instance if $H_G$ is the Hilbert function of a Gorenstein Artinian standard graded $k$-algebra $R/I_G$ then $H_G$ is a Weak Lefschetz sequence.
\par
In~\cite{HMNW}, Proposition 3.5, was shown that $H$ is the Hilbert function of an Artinian  standard graded $k$-algebra having the Weak Lefschetz property iff $H$ is a Weak Lefschetz sequence.
\par
In~\cite{BK}, Example 3.1, was proved that there exist almost complete intersection algebras of codimension $3$ which have not the Weak Lefschetz property, for instance one of this is $R/(x_1^3,x_2^3,x_3^3,x_1x_2x_3)$.
\par
However, in this paper we will prove that all Hilbert functions, which occur for almost complete intersection algebras of codimension $3,$ are Weak Lefschetz sequences.

\section{Main result} 
\markboth{\it Main result}
{\it Main result}

The goal of this section is to show that the Hilbert functions of the almost complete intersection Artinian quotients of $R=k[x_1,x_2,x_3]$ are Weak Lefschetz sequences, i.e. they are unimodal and the positive part of their first differences is an $O$-sequence. 

Let $H=H_Q$ be the Hilbert function of $R/I_Q$ where $I_Q$ is an almost complete intersection ideal of codimension $3.$ Then $I_Q=I_Z:I_G$ for suitable ideals $I_Z \subset I_G$ where $I_Z$ is a complete intersection ideal of codimension $3$ and $I_G$ is a Gorenstein ideal of codimension $3.$ Now, by liaison, we have, for every $i,$ 
$$H_Q(i)=H_Z(\vartheta_Z -3-i)- H_G(\vartheta_Z -3-i)$$ 
which implies
\begin{equation}\label{sim}\Delta H_Q(i)=\Delta H_G(\vartheta_Z -2-i)- \Delta H_Z(\vartheta_Z -2-i).\end{equation}


We start by proving the following general lemma.

\begin{lem}
Let $R=k[x_0,\ldots,x_n]$ be a polynomial ring. Let $I_Z\subseteq I_X\subset R$ be ideals and let $A_Z=R/I_Z$ and $A_X=R/I_X$ be a standard $k$-algebra. Let $l\in R_1$ and let $n\in\zz^+$ such that the multiplication $$l:A_{Z,n-1}\to A_{Z,n}$$ is injective.  Then
 $$\Delta H_{A_Z}(n)\ge\Delta H_{A_X}(n).$$
\end{lem}
\begin{proof}
We set $J_Z=I_Z+(l),$ $J_X=I_X+(l),$ $B_Z=R/J_Z$ and $B_X=R/I_X.$
We have the exact sequence of $k$-vector spaces
 $$0\to A_{Z,n-1}\stackrel{l}{\rightarrow}A_{Z,n}\to B_{Z,n}\to 0,$$
so $H_{B_Z}(n)=\Delta H_{A_Z}(n).$ Since $J_Z\subseteq J_X,$ we have a surjection
 $$B_Z\to B_X\to 0,$$
that implies $H_{B_Z}(n)\ge H_{B_X}(n).$ Now, let $K$ be the kernel of the map $$l:A_X\to A_X.$$
By the exact sequence
 $$0\to K_{n-1}\to A_{X,n-1}\stackrel{l}{\rightarrow}A_{X,n}\to B_{X,n}\to 0,$$
we get $$H_{B_X}(n)=\Delta H_{A_X}(n)+H_K(n-1),$$ therefore
 $$\Delta H_{A_Z}(n)=H_{B_Z}(n)\ge H_{B_X}(n)=\Delta H_{A_X}(n)+H_K(n-1)\ge\Delta H_{A_X}(n).$$ 
\end{proof}

From this lemma we get

\begin{prp}\label{prpzx}
Let $I_X\subset R=k[x_1,x_2,x_3]$ be an Artinian ideal. Let $$f_1,f_2,f_3\in I_X$$ forms such that $(f_1,f_2,f_3)$ is a regular succession. Let $I_Z=(f_1,f_2,f_3).$ We set $A_Z=R/I_Z$ and $A_X=R/I_X.$ Then
 $$\Delta H_{A_Z}(n)\ge\Delta H_{A_X}(n)\,\,\,\,{\text for}\,\,\,\,0\le n\le\vartheta_Z-3-\lambda_Z.$$ 
\end{prp}
\begin{proof}
Since $A_Z$ is a weak Lefschetz Artinian algebra by \cite{HMNW}, there exists a linear form $l\in R_1$ such that the map
 $$l:A_{Z,n-1}\to A_{Z,n}$$
is injective for $0\le n\le\vartheta_Z-3-\lambda_Z.$ So we can apply the previous lemma.
\end{proof}

At this point we set the following notation.

\begin{ntn}
Let $I_G\subset R=k[x_1,x_2,x_3]$ be a Gorenstein Artinian ideal. Let $I_Z\subset I_G$ be an Artinian complete intersection ideal ($I_Z\ne I_G$). Then $\vartheta_G<\vartheta_Z,$ so $\Delta H_Z(\vartheta_Z-2)=-1<\Delta H_G(\vartheta_Z-2)=0.$ We set
 $$\tgz:=\min\{n\in\zz^+_0\mid\Delta H_G(n)>\Delta H_Z(n)\}.$$
\end{ntn}

\begin{prp}\label{lamthe}
Let $I_G\subset R=k[x_1,x_2,x_3]$ be a Gorenstein Artinian ideal. Let $I_Z\subseteq I_G$ be an Artinian complete intersection  ideal. Then
 $$\vartheta_Z-\vartheta_G\ge\lambda_Z-\lambda_G.$$
\end{prp}
\begin{proof}
If $\lambda_Z=\lambda_G$ then $\vartheta_Z-\vartheta_G\ge0=\lambda_Z-\lambda_G;$ so we can suppose $\lambda_Z>\lambda_G.$
We know that $\vartheta_G\ge2\lambda_G+3.$
\par
If $\vartheta_G=2\lambda_G+3,$ then $\vartheta_Z-\vartheta_G\ge2\lambda_Z+3-2\lambda_G-3=2(\lambda_Z-\lambda_G).$
\par
If $\vartheta_G=2\lambda_G+4,$ then $$\vartheta_Z-\vartheta_G\ge2\lambda_Z+3-2\lambda_G-4=
2(\lambda_Z-\lambda_G)-1\ge\lambda_Z-\lambda_G.$$
\par
If $\vartheta_G\ge2\lambda_G+5,$ by Remark~\ref{rmgz}, item $4,$ we have that $\Delta H_G(n)=0$ in at least two integers between $\lambda_G+1$ and $\vartheta_G-3-\lambda_G;$ let $I_Z\in\Gor(a_1,a_2,a_3).$ Then, using Corollary 3.4 in \cite{RZ1}, one sees that $\depth (I_G)_{\le \vartheta_G-3-\lambda_G}=2,$ therefore $a_3\ge\vartheta_G-2-\lambda_G$ so, by Remark~\ref{rmgz}, item $5,$
 $$\vartheta_Z-\vartheta_G=a_1+a_2+a_3-\vartheta_G\ge a_1+a_2-2-\lambda_G\ge\lambda_Z-\lambda_G.$$
\end{proof}

\begin{prp}\label{bvuoto}
Let $I_G\in\Gor\delta,$ $\delta=(d_1,\ldots,d_{2m+1}).$ Let $(d_1,d_{\beta},d_{\gamma})=\mci\delta,$ $d_1\le d_{\beta}\le d_{\gamma}.$
Let us suppose that
 $$B_{\delta}=\{3\le i\le m+1\mid\vartheta_G\le d_i+d_{2m+4-i}\}\ne\emptyset.$$
Then
\begin{itemize}
	\item[1)]$\vartheta_G\le d_{\beta}+d_{\gamma};$
	\item[2)]if $d_1,\ldots,d_{2m+1}$ are the minimal generators degrees allowed by the Hilbert function of $R/I_G,$ then $$\vartheta_G<d_{\beta}+d_{\gamma}.$$
\end{itemize}
\end{prp}
\begin{proof}
\begin{itemize}
\item[1)]If $B_G\ne\emptyset$ then $\beta=\max B_{\delta}$ and $\gamma=2m+4-\min B_{\delta};$ so
 $$2m+4-\gamma\le\beta<2m+4-\beta\le\gamma;$$
consequently we have
 $$d_{2m+4-\gamma}\le d_\beta\le d_{2m+4-\beta}\le d_\gamma.$$
Since $\beta\in B_{\delta}$ we have
 $$\vartheta_G\le d_\beta+d_{2m+4-\beta}\le d_\beta+d_{\gamma}.$$
\item[2)]If $d_1,\ldots,d_{2m+1}$ are the minimal generators degrees allowed by the Hilbert function of $R/I_G,$ then $\vartheta_G\ne d_h+d_k$ for every $h\ne k.$ Since $\beta<\gamma$ we are done.
\end{itemize}
\end{proof}

\begin{prp}\label{riduzb}
Let $I_Q\in\Aci(e_1,e_2,e_3,e_4).$ Then there exists an ideal $I_{Q'}\in\Aci(e_1,e_2,e_3,e_4),$ with $H_{Q'}=H_Q,$ such that $I_{Q'}=I_{Z}:I_{G}$ where $I_{Z}\in\Gor(e_2,e_3,e_4),$ $I_{G}\in\Gor\delta,$ for some $\delta$ with $B_{\delta}=\emptyset,$ $I_{Z}\subset I_{G}$ and $\vartheta_Z-\vartheta_G=e_1.$
\end{prp}
\begin{proof}
$I_Q$ contains a regular sequence of minimal generators of degrees $e_2,e_3,e_4.$ Let $I_{Z_0}$ be the Artinian complete intersection ideal generated by this regular sequence and let 
$I_{G_0}=I_{Z_0}:I_Q.$ Then $I_{G_0}\in\Gor\delta,$ for some $\delta,$ and $\vartheta_{Z_0}-\vartheta_{G_0}=e_1\le e_2$ i.e. $e_2+e_3+e_4-\vartheta_{G_0}\le e_2,$ therefore 
$\vartheta_{G_0}\ge e_3+e_4.$ If $B_{\delta}\ne\emptyset,$ $\mci\delta=(d_1,d_{\beta},d_{\gamma}),$ $d_1\le d_{\beta}\le d_{\gamma},$ then $\vartheta_{G_0}\ge e_3+e_4\ge d_{\beta}+d_{\gamma}$ so
by Proposition~\ref{bvuoto} we get $\vartheta_{G_0}=d_{\beta}+d_{\gamma}.$ So $I_{G_0}$ has not the minimal generators allowed by the Hilbert function of $R/I_{G_0}.$ Now let $I_G\in\Gor\widehat{\delta}$ be another Artinian Gorenstein ideal having the minimal generators allowed by the Hilbert function of $R/I_{G_0}.$ By Proposition~\ref{crescdelta}, item 2, there exists $I_Z\in\Gor(e_2,e_3,e_4)$ such that $I_Z\subset I_G.$ Then, by item $2)$ of Proposition \ref{bvuoto}, $B_{\widehat{\delta}}=\emptyset,$ $I_{Q'}=I_Z:I_G\in\Aci(e_1,e_2,e_3,e_4)$ and $H_{Q'}=H_Q.$ Since $I_Z\in\Gor(e_2,e_3,e_4)$ we will have $e_1=\vartheta_Z-\vartheta_G.$
\end{proof}

Note that, according to Proposition \ref{prpzx} and with the same notation, we have immediately: $\tgz> \vartheta_Z-3-\lambda_Z.$ Now, since $ \vartheta_Z-\lambda_Z>  \vartheta_G-\lambda_G$ (by Proposition \ref{lamthe}), we see also that $\tgz> \vartheta_G-3-\lambda_G.$


\begin{rem}\label{ridg}
Using Proposition~\ref{crescdelta} and Proposition~\ref{riduzb}, in order to produce all the possible Hilbert functions of  almost complete intersection Artinian standard graded $k$-algebras $R/I_Q,$ we can limit us to consider ideals $I_Q=I_Z:I_G$ with  $(I_Z,I_G)\in \Regor(\alpha, \delta),$ $\alpha=(a_1,a_2,a_3),$ such that $B_{\delta}=\emptyset,$ $I_G$ is generated by the minimal generators allowed by the Hilbert function of $R/I_G,$ and $\vartheta_Z-\vartheta_G\le a_1.$
\end{rem}
Now we are ready to prove our main result.

Take $(I_Z,I_G) \in \Regor(\alpha, \delta)$ with $\alpha=(a_1,a_2,a_3))$ and $\delta=(d_1,\ldots,d_{2m+1})$ and  $I_G$ enjoying the properties stated in Remark \ref{ridg} 

\vspace{.2cm}
CLAIM: $\quad \tgz> \vartheta_G-d_2-1.$

\vspace{.2cm}
The proof of the claim is postponed to the next section.

\begin{thm}\label{MT}
Let $H$ be a Hilbert function of an almost complete intersection Artinian standard graded $k$-algebra of codimension $3.$ Then $H$ is a Weak Lefschetz sequence.
\end{thm}
\begin{proof}
Let $R/I_Q$ be a quotient of $R=k[x_1,x_2,x_3]$ with $I_Q$ an almost complete intersection ideal of codimension $3,$ such that $H_Q=H.$ Then $I_Q=I_Z:I_G$ with  $(I_Z,I_G) \in \Regor(\alpha, \delta)$ with $\alpha=(a_1,a_2,a_3),$  $\delta=(d_1,\ldots,d_{2m+1})$ and such that $B_{\delta}=\emptyset,$ $I_G$ is generated by the minimal generators allowed by the Hilbert function of $R/I_G,$ and $\vartheta_Z-\vartheta_G\le a_1$ (see Remark \ref{ridg}). 

Thus, in order to show that $H_Q$ is unimodal we need to find an integer $u$ such that $H_Q(i-1)<H_Q(i)$ for $0< i\le u$ and $H_Q(i-1)\ge H_Q(i)$ for $i> u,$ or, equivalently,
$\Delta H_Q(i)>0$ for $0\le i \le u$ and $\Delta H_Q(i)\le 0$ for $i> u,$ hence, by the equality (\ref{sim}), $\Delta H_G(j)> \Delta H_Z(j)$ for $\vartheta_Z -2-u \le j \le \vartheta_Z -2$ and $\Delta H_G(j)\le \Delta H_Z(j)$ for $0 \le j < \vartheta_Z -2-u.$

Define $u:=\vartheta_Z-2-\tgz;$ for every $ 0\le i < \tgz$, by definition of $\tgz,$  we have $\Delta H_G(i)\le \Delta H_Z(i),$ i.e. $\Delta H_G(i)\le \Delta H_Z(i)$ for $0 \le i < \vartheta_Z -2-u.$ Now, we want to show that for every $i,$  $\vartheta_Z -2-u \le i \le \vartheta_Z -2$ we have $\Delta H_G(i)> \Delta H_Z(i).$ If not, let $j=\min\{i \ | \ \vartheta_Z -2-u < i \le \vartheta_Z -2, \ \Delta H_G(i)\le \Delta H_Z(i)\},$ then  $\Delta H_G(j-1)> \Delta H_Z(j-1)\}.$ Now, since $j>\vartheta_Z -2-u=\tgz,$ by the previous claim, $j >\vartheta_G-d_2-1,$ thus $\Delta^2H_G(j)\ge 0.$ If $j\le \vartheta_Z-a_1-1$ then $\Delta^2H_Z(j)\le 0.$ From this, we get 
$$\Delta H_G(j)=\Delta H_G(j-1)+\Delta^2 H_G(j)> \Delta H_Z(j-1)+\Delta^2 H_Z(j)=\Delta H_Z(j)$$ 
a contradiction. If  $ \vartheta_Z-a_1-1 <j \le \vartheta_G -2$ then $\Delta^2H_G(j)=\Delta^2H_Z(j)=1,$ so again
$$\Delta H_G(j)=\Delta H_G(j-1)+1> \Delta H_Z(j-1)+1=\Delta H_Z(j).$$
Finally, if $\vartheta_G-2 <j \le \vartheta_Z -2,$ trivially  $0=\Delta H_G(j)>\Delta H_Z(j).$ In this way we can conclude that $H$
is unimodal. 

In order to conclude the proof we need to show that the sequence $(\Delta H)^+$  is an $O$-sequence. In our situation, using Proposition \ref{riduzb}, it will be enough to show that, for $i\ge \vartheta_Z - \vartheta_G,$ $(\Delta H)^+(i-1)\ge (\Delta H)^+(i).$ This is equivalent to prove that, for $\tgz+1\le j\le \vartheta_G-1,$ $\Delta^2 H_G(j)\ge \Delta^2 H_Z(j).$ Since $\Delta^2 H_G(j)=1$ for $j > \vartheta_G-1-d_1$ the inequality trivially holds in this range. Now, take  $\tgz+1\le j\le \vartheta_G-1-d_1.$ Since  $\vartheta_G < d_1+a_2+a_3,$
as $\mci(d_1, \ldots, d_{2m+1})\le (d_1,a_2,a_3),$ we have that $j \le \vartheta_Z-a_1-1,$ therefore $\Delta^2 H_G(j)=0$ and $\Delta^2 H_Z(j)\le 0$ so we get again the required inequality.\end{proof}

\section{Proof of the claim} 
\markboth{\it Proof of the claim}
{\it Proof of the claim}
Let $I_Z\subset I_G$ be ideals such that $I_Z\in\Gor\alpha,$ $\alpha=(a_1,a_2,a_3),$ $I_G\in\Gor\delta,$ $\delta=(d_1,\ldots,d_{2m+1})$ and $I_G$ satisfying the conditions of Remark \ref{ridg}.
\par
In this section we will prove the claim stated in the Section $2$, i.e.
 $$\tgz> \vartheta_G-d_2-1,$$
which means to prove that 
 $$\Delta H_G(i)\le\Delta H_Z(i),\,\text{ for }0\le i\le\vartheta_G-d_2-1.$$

At first we prove the claim when $I_G$ is a complete intersection  ideal.

\begin{prp}\label{gcoint}
Let $(I_Z,I_G)\in\Regor(\alpha,\delta),$ with $\alpha=(a_1,a_2,a_3)$ and $\delta=(d_1,d_2,d_3).$ Then
 $$\Delta H_G(i)\le\Delta H_Z(i),\,\text{ for }0\le i\le\vartheta_G-d_2-1.$$
\end{prp}
\begin{proof} 
By hypotheses $d_i\le a_i,$ for $1\le i\le 3.$
\par
If $d_3\le d_1+d_2-1$ then 
$$\Delta H_G(i)=\left\{
\begin{array}{cl}
	i+1 & \text{ for } 0\le i\le d_1-1\\
	d_1   & \text{ for } d_1-1\le i\le d_2-1 \\
	-i+d_1+d_2-1   & \text{ for }  d_2-1\le i\le d_3-1 \\
	-2i+d_1+d_2+d_3-2 & \text{ for }  d_3-1\le i\le d_1+d_2-1 \\
	-i+d_3-1 & \text{ for }  d_1+d_2-1\le i\le d_1+d_3-1 \\
	-d_1 & \text{ for }  d_1+d_3-1\le i\le d_2+d_3-1 \\
	i-d_1-d_2-d_3+1 & \text{ for }  d_2+d_3-1\le i\le d_1+d_2+d_3-1 \\
	0 & \text{ for }  i\ge d_1+d_2+d_3-1
\end{array}\right.;$$
if $d_1+d_2-1\le d_3$ then
$$\Delta H_G(i)=\left\{
\begin{array}{cl}
	i+1 & \text{ for } 0\le i\le d_1-1\\
	d_1   & \text{ for } d_1-1\le i\le d_2-1 \\
	-i+d_1+d_2-1   & \text{ for }  d_2-1\le i\le d_1+d_2-1 \\
	0 & \text{ for }  d_1+d_2-1\le i\le d_3-1 \\
	-i+d_3-1 & \text{ for }  d_3-1\le i\le d_1+d_3-1 \\
	-d_1 & \text{ for }  d_1+d_3-1\le i\le d_2+d_3-1 \\
	i-d_1-d_2-d_3+1 & \text{ for }  d_2+d_3-1\le i\le d_1+d_2+d_3-1 \\
	0 & \text{ for }  i\ge d_1+d_2+d_3-1
\end{array}\right.;$$
and analogously for $\Delta H_Z.$
Therefore the assertion follows after straightforward computations.
\end{proof}

\begin{lem}\label{ridmin}
Let $(I_Z,I_G)\in\Regor(\alpha,\delta),$ with $\alpha=(a_1,a_2,a_3)$ and $\delta=(d_1,\ldots,d_{2m+1}),$ $m\ge 2,$ $B_{\delta}=\emptyset.$ Let $\alpha_0=\min\delta$ and let $I_{Z_0}\subset I_G,$ $I_{Z_0}\in\Gor\alpha_0.$ Let us suppose that
 $$\Delta H_G(i)\le\Delta H_{Z_0}(i),\,\text{ for }0\le i\le\vartheta_G-d_2-1.$$
Then
 $$\Delta H_G(i)\le\Delta H_{Z}(i),\,\text{ for }0\le i\le\vartheta_G-d_2-1.$$
\end{lem}
\begin{proof} 
$B_{\delta}=\emptyset$ implies that $\alpha_0=(d_1,d_2,d_{\gamma}),$ for some $\gamma \ge 3.$ Since $\alpha\ge\alpha_0,$ $\Regor(\alpha,\alpha_0)\ne\emptyset.$ Let $(I_W,I_{W_0})\in\Regor(\alpha,\alpha_0).$ By Proposition \ref{gcoint} we have
 $$\Delta H_{W_0}(i)\le\Delta H_W(i),\,\text{ for }0\le i\le\vartheta_{W_0}-d_2-1$$ i.e.
 $$\Delta H_{Z_0}(i)\le\Delta H_Z(i),\,\text{ for }0\le i\le\vartheta_{Z_0}-d_2-1.$$
But $\vartheta_G-d_2-1\le\vartheta_{Z_0}-d_2-1,$ so we have
 $$\Delta H_G(i)\le\Delta H_{Z}(i),\,\text{ for }0\le i\le\vartheta_G-d_2-1.$$
\end{proof}

\begin{prp}\label{cvuoto}
Let $I_Z\subset I_G$ be ideals such that $I_Z\in\Gor\alpha,$ $\alpha=(a_1,a_2,a_3),$ $I_G\in\Gor\delta,$ $\delta=(d_1,\ldots,d_{2m+1}),$ $m\ge 2,$ $I_G$ satisfying the conditions of Remark \ref{ridg} and $C_{\delta}=\emptyset.$
Then $\Delta H_G(i)\le\Delta H_Z(i),$\, for $0\le i\le\vartheta_G-d_2-1.$\end{prp}
\begin{proof} Since $B_{\delta}=C_{\delta}=\emptyset,$ $\mci\delta=(d_1,d_2,d_3).$ By Lemma \ref{ridmin}, it is enough to prove the inequality for $\alpha=(d_1,d_2,d_3).$ Because of that $\Delta H_G(i)=\Delta H_{Z}(i)$ for $0\le i\le d_3-1.$ Now, observe that, since  
$B_{\delta}=C_{\delta}=\emptyset,$ $\Delta^2H_G(i) \le -2$ for $d_3 \le i \le \vartheta_G -d_3 -1$ (see Proposition 3.11 in~\cite{RZ2}). Reminding that  $\Delta^2H_Z(i) \ge -2$ for all $i,$ we  get immediately that the inequality holds for $0\le i\le \vartheta_G -d_3 -1.$ Now let $p:=\Delta H_G(\vartheta_G-d_3-1)=\Delta H_Z(\vartheta_Z-d_3-1).$ $\Delta H_G(i)<p$ for 
$\vartheta_G-d_3\le i\le\vartheta_G-d_2-1$ and $p\le \Delta H_Z(i)$ for $i\le\vartheta_Z-d_3-1;$ so the assertion follows if $\vartheta_G-d_2-1\le\vartheta_Z-d_3-1.$ Otherwise, if $\vartheta_G-d_2-1>\vartheta_Z-d_3-1,$ we have
$\Delta H_G(i)<p\le\Delta H_Z(i)$ for $\vartheta_G-d_3\le i\le\vartheta_Z-d_3-1.$ Finally, if $i$ is any integer such that $\vartheta_Z-d_3-1\le i\le\vartheta_G-d_2-1,$ we have $\Delta H_G(i)=\vartheta_G-(i+1+d_3-p)$ and $\Delta H_Z(i)=\vartheta_Z-(i+1+d_3-p),$ so, again, $\Delta H_G(i)<\Delta H_Z(i).$
\end{proof}

\begin{prp}\label{cnonvuoto}
Let $I_Z\subset I_G$ be ideals such that $I_Z\in\Gor\alpha,$ $\alpha=(a_1,a_2,a_3),$ $I_G\in\Gor\delta,$ $\delta=(d_1,\ldots,d_{2m+1}),$ $m\ge 2,$ $I_G$ satisfying the conditions of Remark \ref{ridg} and $C_{\delta}\ne\emptyset.$
Then $\Delta H_G(i)\le\Delta H_Z(i),$\, for $0\le i\le\vartheta_G-d_2-1.$
\end{prp}
\begin{proof} Since $B_{\delta}=\emptyset$ and $C_{\delta}\ne\emptyset,$ $\mci\delta=(d_1,d_2,d_{\gamma}),$ with $\gamma\ge 4.$ By Lemma \ref{ridmin}, it is enough to prove the inequality for $\alpha=(d_1,d_2,d_{\gamma}).$ Because of that $\Delta H_G(i)=\Delta H_{Z}(i)$ for $0\le i\le d_2-1.$ Now, observe that, since  
$B_{\delta}=\emptyset$ and $C_{\delta}\ne\emptyset,$ $\Delta^2H_G(i) \le -1$ for $d_2\le i\le d_{\gamma}-1$ (see Proposition 3.11 in~\cite{RZ2}), and $\Delta^2H_Z(i)=-1$ in the same range; so the inequality holds for $0\le i\le d_{\gamma}-1.$ 
\par
Let $p:=\Delta H_G(d_{\gamma}-1);$ then $\Delta H_G(i)\le -2i+p+2d_{\gamma}-2$ for 
$d_{\gamma}-1\le i\le\vartheta_G-d_{\gamma}-1$ and $\Delta H_G(i)\le -i+\vartheta_G-p-d_{\gamma}-1$ for 
$\vartheta_G-d_{\gamma}-1\le i\le\vartheta_G-d_{2}-1.$ Moreover 
$-p=\Delta H_G(\vartheta_G-d_{\gamma}-1)\le -2(\vartheta_G-d_{\gamma}-1)+p+2d_{\gamma}-2$ i.e. 
$\vartheta_G\le p+2d_{\gamma}.$ Since $\Delta^2H_G(i)\le -2$ for $d_{\gamma}\le i\le\vartheta_G-d_{\gamma}-1$ we deduce that
 $\Delta H_G(i)\le\Delta H_Z(i),$\, for $0\le i\le\vartheta_G-d_{\gamma}-1.$ It remains to show the inequality for $\vartheta_G-d_{\gamma}-1\le i\le\vartheta_G-d_{2}-1.$ 
\par
If $\vartheta_Z-d_{\gamma}-1\le\vartheta_G-d_{2}-1,$ take $i$ such that $\vartheta_G-d_{\gamma}-1\le i\le\vartheta_Z-d_{\gamma}-1;$ then 
 $$\Delta H_G(i)\le -i+\vartheta_G-p-d_{\gamma}-1\le-2i+\vartheta_Z-2\le\Delta H_Z(i),$$
since $i\le\vartheta_Z-d_{\gamma}-1$ and $\vartheta_G\le p+2d_{\gamma}.$ Now take $i$ such that $\vartheta_Z-d_{\gamma}-1\le i\le\vartheta_G-d_{2}-1;$ then
 $$\Delta H_G(i)\le -i+\vartheta_G-p-d_{\gamma}-1\le-i+d_{\gamma}-1=\Delta H_Z(i),$$
since $\vartheta_G\le p+2d_{\gamma}.$
\par
If $\vartheta_G-d_{2}-1<\vartheta_Z-d_{\gamma}-1,$ take $i$ such that $\vartheta_G-d_{\gamma}-1\le i\le\vartheta_G-d_{2}-1;$ then
 $$\Delta H_G(i)\le -i+\vartheta_G-p-d_{\gamma}-1\le-2i+\vartheta_Z-2\le\Delta H_Z(i),$$
since, as above, $i\le\vartheta_Z-d_{\gamma}-1$ and $\vartheta_G\le p+2d_{\gamma}.$
\end{proof}

\vspace{1cm}
{\f
{\sc (A. Ragusa) Dip. di Matematica e Informatica, Universit\`a di Catania,\\
                  Viale A. Doria 6, 95125 Catania, Italy}\par
{\it E-mail address: }{\tt ragusa@dmi.unict.it} \par
{\it Fax number: }{\f +39095330094} \par
\vspace{.3cm}
{\sc (G. Zappal\`a) Dip. di Matematica e Informatica, Universit\`a di Catania,\\
                  Viale A. Doria 6, 95125 Catania, Italy}\par
{\it E-mail address: }{\tt zappalag@dmi.unict.it} \par
{\it Fax number: }{\f +39095330094}
}

\end{document}